\documentclass[aap,MSNbibl,seceqn,nameyear,dvips]{arximspdf}
\usepackage{mathbh}

\doi{10.1214/12-AAP908} 
\volume{23}
\issue{6}
\pubyear{2013}
\firstpage{2472}
\lastpage{2499}

\makeatletter

\newcommand{\cal}{\mathcal}

\def\Om{\Omega}
\def\ua{\uparrow}

\def\bR{\mathbb{R}}

\def\Om{\Omega}

\def\cA{{\cal A}}
\def\cB{{\cal B}}
\def\cF{{\cal F}}
\def\cG{{\cal G}}
\def\cM{{\cal M}}

\def\lra{\longrightarrow}
\def\ua{\uparrow}

\def\bR{\mathbb R}
\def\cF{\mathcal F}
\def\bP{\mathbb P}
\def\bE{\mathbb E}
\def\bN{\mathbb N}
\def\var{\operatorname{var}}

\newtheorem{theorem}{Theorem}[section]
\newtheorem{proposition}[theorem]{Proposition}
\newtheorem{lemma}[theorem]{Lemma}

\newproclaim{example}[theorem]{Example}
\newproclaim{remark}[theorem]{Remark}

\def\Ind#1{\mathbh1_{_{#1}}}

\makeatother

\begin{document}
\begin{frontmatter}

\title{A control problem with fuel constraint and Dawson--Watanabe
superprocesses}
\runtitle{A control problem and superprocesses}

\begin{aug}
\author[A]{\fnms{Alexander} \snm{Schied}\corref{}\thanksref{t1}\ead[label=e1]{schied@uni-mannheim.de}}
\runauthor{A. Schied}
\affiliation{University of Mannheim}
\address[A]{Department of Mathematics\\
University of Mannheim\\
A5,6\\
68131 Mannheim\\
Germany\\
\printead{e1}} 
\end{aug}

\thankstext{t1}{Supported by the Deutsche Forschungsgemeinschaft (DFG).}

\received{\smonth{7} \syear{2012}}
\revised{\smonth{10} \syear{2012}}

%
\begin{abstract}
We solve a class of control problems with fuel constraint by means of
the log-Laplace transforms of $J$-functionals of Dawson--Watanabe
superprocesses. This solution is related to the superprocess solution
of quasilinear parabolic PDEs with singular terminal condition. For the
probabilistic verification proof, we develop sharp bounds on the
blow-up behavior of log-Laplace functionals of $J$-functionals, which
might be of independent interest.
\end{abstract}

%

\end{frontmatter}

\section{Introduction}

One of the most exciting aspects of Dawson--Watanabe superprocesses is
their connection to quasilinear partial differential equations (PDEs)
with singular boundary condition. This connection was pioneered by
Dynkin (\citeyear{DynkinElliptic,DynkinParabolic}); see also \citet{DynkinPDEBook}
for more recent developments and related literature. Similar
quasilinear PDEs also appear in the Hamilton--Jacobi--Bellman (HJB)
formulation of stochastic control problems with terminal state
constraint, and so it is natural to ask whether these control problems
possess solutions in terms of superprocesses. Establishing such a
direct connection is the main goal of this paper.

The connection we find is direct insofar as it avoids the use of HJB
equations and instead uses a probabilistic verification argument based
solely on the log-Laplace equation for a certain $J$-functional of a
superprocess. While the standard verification argument relies on the
existence of smooth solutions to the HJB equation, whose existence is
often very difficult to establish, the ``mild solutions'' provided by
the log-Laplace functionals of superprocesses are ideally suited for
carrying out the verification argument. They are also superior to the
commonly used viscosity solutions, because the latter do not go well
along with It\^o calculus due to their possible lack of smoothness.

The problem we will consider here is the minimization of the functional
%
\begin{equation}
\label{introproblem} E_{0,z} \biggl[ \int_0^T\bigl|
\dot x(t)\bigr|^p\eta(Z_t)\,dt+\int_{[0,T]}\bigl|x(t)\bigr|^p
A(dt) \biggr]
\end{equation}
over adapted and absolutely continuous strategies $x(t)$ satisfying the
constraints $x(0)=x_0$ and $x(T)=0$. We assume here that $p\in
[2,\infty)$, $\eta$ is a strictly positive function and $A$ is a
nonnegative additive functional of the (time-inhomogeneous) Markov
process $Z$ with $Z_0=z$ $P_{0,z}$-a.s. This control problem is closely
related to the monotone follower problems with fuel constraint that
were introduced by \citet{BenesSheppWitsenhausen} and further
developed, for example, by \citet{Karatzas}. Also, as we will explain
in Section \ref{financesection}, problems of this type have recently
appeared in the context of mathematical finance.
In the next section we will give some heuristic arguments that explain
the connection of this problem with quasilinear PDEs with singular
terminal condition that are related to superprocesses.

\subsection{The connection between the control problem and superprocesses}

Let us assume that $A(du)=a(Z_u)\,du$ for some function $a\ge0$ and
define the value function of our problem as
\[
V(t,z,x_0):=\inf_{x(\cdot)}E_{t,z} \biggl[
\int_t^T\bigl|\dot x(u)\bigr|^p\eta
(Z_u)\,du+\int_t^T\bigl|x(u)\bigr|^pa(Z_u)
\,du \biggr],
\]
where the infimum is taken over the class of all absolutely continuous
and adapted strategies $x(\cdot)$ satisfying the constraints
$x(t)=x_0$ and $x(T)=0$. As usual $P_{t,z}$
denotes the probability measure under which the Markov process $Z$
starts at $z$ at time $t$.
When $L_t$ is the infinitesimal generator of $Z$, standard arguments
from optimal control suggest that $V$ should satisfy the following HJB equation:
%
\begin{eqnarray}
\label{VHJBeqn}
&&V_t(t,z,x_0)+\inf_{\xi}
\bigl\{\eta(z)|\xi|^p+V_{x_0}(t,z,x_0)\xi
\bigr\}\nonumber\\[-8pt]\\[-8pt]
&&\qquad{}+a(z)|x_0|^p+L_tV(t,z,x_0)=0\nonumber
\end{eqnarray}
with singular terminal condition
%
\begin{equation}
\label{Vterminalcondition} V(T,z,x_0)=\cases{0, &\quad if
$x_0=0$,
\cr
+\infty, &\quad otherwise.}
\end{equation}
Note that the singularity in this terminal condition is required by the
fuel constraint $x(T)=0$.

To see how this PDE is related to superprocesses, we consider the case
$x_0\ge0$ and make the ansatz $V(t,z,x_0)=x_0^pv(t,z)$ for some
function $v$. Plugging this ansatz into (\ref{VHJBeqn}), minimizing
over $\xi$, dividing by $x_0^p$ and using (\ref{Vterminalcondition})
yields the equation
%
\begin{eqnarray}
\label{DynkinPDE} v_t-\frac1{\beta\eta^{\beta}}v^{1+\beta}+a+L_tv&=&0,
\nonumber\\[-8pt]\\[-8pt]
v(T,z)&=&+\infty,
\nonumber
\end{eqnarray}
where
$\beta=\frac1{p-1}$. This is just the type of PDE solved in \citet
{DynkinParabolic} by means of superprocesses.

The minimizing $\xi$ in (\ref{VHJBeqn}) is given by $\xi=-x_0
v^\beta/\eta^\beta$, which suggests that the minimizing strategy
$x^*$ for the cost functional (\ref{introproblem}) is given in
feedback form as the solution of the ordinary differential equation
\[
\dot x(u)=-\frac{x(u)v(u,Z_u)^\beta}{\eta(Z_u)^\beta},
\]
that is,
\[
x^*(t)=x_0\exp \biggl(-\int_0^t
\frac{v(s,Z_s)^\beta}{\eta(Z_s)^\beta
}\,ds \biggr).
\]

As we will see later on, these heuristic computations will give the
correct results. There are some difficulties, however, which need to be
overcome to turn these heuristics into a full proof. For instance, we
must make sure that the strategy $x^*$ satisfies the fuel constraint
$x^*(T)=0$, which requires us to find a \textit{lower} bound on
$v(s,Z_s)$ when it approaches the singularity so as to ensure that
$\int_0^t\frac{v(s,Z_s)^\beta}{\eta(Z_s)^\beta}\,ds$ diverges when
$t\ua T$. On the other hand, we must also make sure that $x^*$ has
finite cost (\ref{introproblem}). To this end, we will need a sharp
\textit{upper} bound on $v(s,Z_s)$ for $s$ close to $T$. These bounds
are derived in Section \ref{auxSection} by extending existing bounds
from \citet{SchiedDisspaper} to $J$-functionals and to a generalized
class of superprocesses with nonhomogeneous branching parameters. These
bounds may be of independent interest. For instance, by means of the
results in \citet{DynkinParabolic} they translate into sharp bounds for
solutions $v(t,z)$ of the singular Cauchy problem (\ref{DynkinPDE})
when $t$ approaches the singularity $T$.

\subsection{Financial motivation of the control problem}\label
{financesection}

Let $(S_t)$ be a square-integrable martingale, which will be
interpreted as the price process of an asset. In the linear
Almgren--Chriss market impact model, a strategy $x(\cdot)$ as
described above is interpreted as the policy of an investor who wishes
to liquidate $x_0$ shares of the asset throughout the time interval
$[0,T]$. This liquidation creates price impact so that the investor
trades at price $S^x_t=S_t+\gamma(x_t-x_0)+\eta_t\dot x_t$, where
$\gamma$ is a constant and the process $\eta$ describes the intraday
liquidity fluctuations; see \citet{AlmgrenChriss2}, \citet
{AlmgrenSIFIN} and the survey \citet{GatheralSchiedSurvey} for details.
The liquidation costs arising from the strategy $x(\cdot)$ are then
given by
\[
C(x)=\frac\gamma2x_0^2+\int_0^T
\eta_t\dot x(t)^2\,dt-\int_0^Tx(t)
\,dS_t.
\]
The problem considered by practitioners is the minimization of the
following mean-variance functional of the costs:
%
\begin{equation}
E\bigl[ C(x) \bigr]+\lambda\var\bigl(C(x)\bigr),
\end{equation}
over all absolutely continuous policies satisfying the constraints
$x(0)=x_0$ and $x(T)=0$. This is a straightforward exercise when $\eta
$ is constant and strategies are deterministic, but not so easy when
strategies are adapted. The reason for this difficulty is the time
inconsistency of the mean-variance functional, which precludes the use
of control techniques; see, for example, \citet{Tseetal} and the
references therein.

As a way out, one can use an infinitesimal re-optimization process as
in Section~6.4 of \citet{Schoeneborn} or other, more generally
available arguments such as those in \citet{EkelandLazrak} or \citet
{Bjoerk} to replace the original, time-inconsistent problem by a
time-consistent approximation. At least when $\eta$ is deterministic,
this process leads to the problem of minimizing the functional
%
\begin{equation}
\label{AlmgrenCostEq} E \biggl[ \int_0^T
\eta_t\dot x(t)^2\,dt+\lambda\int_{[0,T]}x(t)^2
\,d[ S]_t \biggr];
\end{equation}
see also \citet{AlmgrenSIFIN}, \citet{Forsythetal} and \citet{Tseetal}
for other motivations of this problem and further studies with
applications in mind. In particular, \citet{AlmgrenSIFIN} proposes to
study the cost functional (\ref{AlmgrenCostEq}) also for
nondeterministic $\eta$.
When $S_t$ and $\eta_t$ are functions of an underlying Markov process
$Z_t$, as it is the case for almost all probabilistic models of asset
price processes, we see that this problem is precisely of the form
(\ref{introproblem}). It is now also clear that using a general
additive functional $A$ and not just a functional of the form
$A(dt)=a(Z_t)\,dt$
in the formulation of (\ref{introproblem}) is suggested by applications
and not just an artificial generalization of our problem.

\subsection{Plan of the paper}

In Section \ref{Assumptionssection} we introduce the basic setup of the
paper. In Section \ref{Statementsection} we first state the solution of
problem (\ref{introproblem}) for $\eta\equiv1$. Theorem
\ref{homogeneousthm}, the corresponding result, is actually a special
case of our main result, Theorem \ref{mainthm}, but unlike Theorem
\ref{mainthm} it only involves classical superprocesses, as constructed
in \citet{DynkinHistorical}, while the case of nonconstant $\eta$
requires an extended class of superprocesses. In Theorem \ref
{relaxedthm} we consider a variant of problem (\ref{introproblem}) in
which the fuel constraint $x(T)=0$ is relaxed and replaced by a penalty
term.

In Section \ref{auxSection} we collect some auxiliary results on
superprocesses and their Laplace functionals, some of which may be of
independent interest. In Section~\ref{maxprinzSection} we present a
probabilistic version of the parabolic maximum principle for the
log-Laplace equations associated with $J$-functionals. This result will
be needed for comparing log-Laplace equations for superprocesses with
inhomogeneous branching characteristics to the case with homogeneous
branching. In Section \ref{HomogeneousLaplaceestimatesSection} we
derive estimates for the Laplace functionals of $J$-functionals of
superprocesses with homogeneous branching characteristics by extending bounds
obtained in \citet{SchiedDisspaper} to the case of $J$-functionals. In
Section \ref{h-trafosection}, these bounds are then extended as well
to a generalized class of superprocesses with nonhomogeneous branching
parameters. This extension is needed for Theorem \ref{mainthm}. The
generalized class of superprocesses is constructed by means of an
``$h$-transform'' for superprocesses that was introduced
independently by \citet{EnglaenderPinsky} and \citet{SchiedDiplompaper}.

The proofs of our main results are given in Section \ref{Proofssection}.

\section{Setup, preliminaries, and main results}

\subsection{Assumptions and preliminaries}\label{Assumptionssection}

\subsubsection{The Markov process $Z$}

We will assume henceforth that the Markov process $Z=(Z_t,\cF
(I),P_{r,z})$ is a time-in\-homo\-geneous right process with sample
space $(\Om,\cF)$ and state space $({S},\cB)$ in the sense of
Section 2.2 in \citet{DynkinBranchingBook}. Here, ${S}$ is a metrizable
Luzin space with Borel field $\cB$. The $\sigma$-algebra $\cF
(I)\subset\cF$ contains events observable during the time interval
$I\subset[0,\infty)$. For every $r\ge0$ and probability measure $\mu
$ on $(S,\cB)$, we thus get a filtered probability space $(\Om,(\cF
[r,t])_{t\ge r},P_{r,\mu})$. Here, $P_{r,\mu}$ denotes as usual the
probability measure under which $Z$ starts at time $r\ge0$ with
initial distribution $\mu$. The fact that $Z$ is right essentially
means that $t\mapsto Z_t(\omega)$ is c\`adl\`ag for each $\omega$ and
that, for $r<t$, measurable $f\dvtx {S}\to\bR_+$, and probability measures
$\mu$ on $(S,\cB)$, the process $s\mapsto E_{s,Z_s}[ f(Z_t) ]$ is
$P_{r,\mu}$-a.s. right-continuous on $[r,t)$. For all further details
we refer to Sections 2.2.1 and 2.2.3 in \citet{DynkinBranchingBook}.

\subsubsection{\texorpdfstring{The $(Z,K,\psi)$-superprocess}
{The (Z, K, psi)-superprocess}}

For $t\ge0$, $z\in{S}$, $\xi\ge0$, $\beta\in(0,1]$ and bounded,
measurable and positive $a\dvtx \bR_+\times{S}\to\bR_+$ let
%
\begin{equation}
\label{psiclass} \psi(t,z,\xi)=a(t,z)\xi^{1+\beta}.
\end{equation}
Let moreover $K$ be a continuous and nonnegative additive functional of
$Z$ such that
%
\begin{equation}
\label{Kcondition} \sup_\omega K[0,t](\omega)<\infty\qquad\mbox{for
all $t\ge0$.}
\end{equation}
By Theorem 1.1 from \citet{DynkinHistorical} one
can then\setcounter{footnote}{1}\footnote{Actually, a larger class of functions $\psi$ is possible
Dynkin (\citeyear{DynkinElliptic,DynkinBranchingBook}), but here we will only need the class specified
in (\ref{psiclass}).} construct the \textit{superprocess with
parameters $(Z,K,\psi)$}. It is a time-inhomogeneous Markov process
$X=(X_t,\cG(I),\bP_{r,\mu})$ with state space $\cM$, the space of
all nonnegative finite Borel measures on $({S},\cB)$. Its transition
probabilities are determined as follows by the Laplace functionals of
$X$. For any measurable function $f\dvtx {S}\to\bR_+$ and $\nu\in\cM$,
we write $\langle f,\nu\rangle$ shorthand for the integral $\int f\,d\nu
$. Then
\[
\bE_{r,\mu}\bigl[ e^{-\langle f,X_T\rangle} \bigr]=e^{-\langle v(r,\cdot
),\mu\rangle},
\]
where $v$ solves the integral equation
%
\begin{eqnarray}
\label{vlog-Laplaceequation} v(s,z)&=&E_{s,z}\bigl[ f(Z_T)
\bigr]\nonumber\\[-8pt]\\[-8pt]
&&{}-E_{s,z} \biggl[ \int_s^T\psi
\bigl(t,Z_t,v(t,Z_t)\bigr) K(dt) \biggr],\qquad r\le s\le T.\nonumber
\end{eqnarray}
When $f$ is bounded, then $v$ is the unique nonnegative solution of
(\ref{vlog-Laplaceequation}). For $T<r$ we make the convention that
$X_T=0$ $\bP_{r,\mu}$-a.s.

Later on, we will need to allow for unbounded additive functionals $K$
and thus have to extend this class of superprocesses; see
Proposition \ref{h-transformexistProp}.

It might be interesting to note that superprocesses provide an
infinite-dimension\-al example of an affine process, a class of processes
which has recently received considerable attention in mathematical
finance; see \citet{DuffieFilipovicSchachermayer}.

\subsubsection{A class of additive functionals of $Z$}

Following \citet{DynkinHistorical}, a~nonnegative additive functional
$A$ of $Z$ belongs to $\cA_{(1)}$ if there exists a finite set $\{
t_1,\ldots, t_n\}\subset\bR_+$ with $t_1<\cdots<t_n$ and bounded
measurable functions $f_i\ge0$, $i=1,\ldots, n$, such that
%
\begin{equation}
\label{ADelta} A[s,u]=\sum_{s\le t_i\le u}f_i(Z_{t_i}).
\end{equation}
Next, $\cA_{(2)}$ denotes the class of all nonnegative additive
functionals $A$ for which there exists a sequence $(A_n)$ in $\cA
_{(1)}$ with the following three properties:  $A_n[r,\infty)\to
A[r,\infty)$ $P_{r,z}$-a.s. for all pairs $(r,z)$; there exists $T>0$
such that $A_n[T,\infty)=0$ for all $n$; and $\sup_{\omega,
n}A_n[0,T](\omega)<\infty$. Finally, $\cA$ consists of all
nonnegative additive functionals $A$ for which there exists a sequence
$(A_n)$ in $\cA_{(2)}$ such that $A_n(B)\nearrow A(B)$ for all
measurable sets $B\subset\bR_+$.
For $q\ge1$ and $T>0$, we furthermore introduce the class
\begin{eqnarray*}
\cA_T^q&:=& \bigl\{A\in\cA| E_{r,z}\bigl[
A[r,T]^q \bigr]<\infty\mbox{ and }\\
&&\hspace*{4.7pt} A(T,\infty)=0\mbox{
$P_{r,z}$-a.s. for all $(r,z)\in[0,T]\times S$} \bigr\}.
\end{eqnarray*}

\begin{remark}[(Quadratic variation and path processes)]
Suppose that $Y_t$ is a semimartingale of the form $Y_t=\phi(Z_t)$,
where $\phi\dvtx {S}\to\bR$ is a measurable function. Then the quadratic
variation of $Y$ gives rise to the nonnegative additive functional
%
\begin{equation}
\label{Aquadrvar}A(dt):=d[Y]_t
\end{equation}
of $Z$. But in general it is not obvious whether $A$ can be
approximated by additive functionals of the form (\ref{ADelta})
unless, for example, $Y$ is an It\^o process of the form $Y_t=Y_0+\int_0^t\sigma(s,Z_s)\,dW_s+\int_0^tb(s,Z_s)\,ds$ and so $d[Y]_t=\sigma
(t,Z_t)^2\,dt$. Nevertheless, when we are interested primarily in $Y$
and its quadratic variation as in the financial context of Section \ref
{financesection}, then we may always assume that $Z$ is the
\textit{path process} for $Y$ [called \textit{historical process} in \citet
{DawsonPerkins} or \citet{Perkins}]; that is, at each point $t$ in
time, $Z_t$ is equal to the sample path of the entire history
$(Y_{s\wedge t})_{s\ge0}$. Then the dynamics of $Z$ will automatically
be Markovian, and $Z$ will be a right process under mild assumptions.
In this case, we can let
\[
f_i(Z_{t_i}):= \bigl(Z_{t_i}({t_i})
-Z_{t_i}({t_{i-1}}) \bigr)^2\wedge K=
(Y_{t_i}-Y_{t_{i-1}} )^2\wedge K
\]
in (\ref{ADelta}),
where $K>0$ and $t_0:=0$. The corresponding additive functional belongs
to $\cA_{(1)}$ and can be used to approximate the quadratic variation
$[Y]$, so that under mild conditions (\ref{Aquadrvar}) has a version
in $\cA$.
\end{remark}

\subsubsection{$J$-functionals associated with an additive functional}

$J$-functionals are functionals of the $(Z,K,\psi)$-superprocess $X$
associated with additive functionals $A\in\cA$. They were introduced
in \citet{DynkinHistorical} as follows. Suppose that $A\in\cA_{(1)}$
is given by (\ref{ADelta}).
Then the corresponding $J$-functional is defined as
%
\begin{equation}
\label{IAdef} J_A:=\sum_{i=1}^n
\langle f_{t_i},X_{t_i}\rangle
\end{equation}
(recall the convention that $X_t=0$ $\bP_{r,\mu}$-a.s. for $t<r$). For
more general additive functionals $A$, the corresponding
$J$-functionals can then be defined by a limiting procedure. By Theorem
1.2 in \citet{DynkinHistorical}, one has for $A\in\cA^1_T$,
%
\begin{equation}
\label{IALaplace} \bE_{r,\mu}\bigl[ e^{-J_A}
\bigr]=e^{-\langle v(r,\cdot),\mu\rangle},
\end{equation}
where $v(r,z)$ solves
%
\begin{equation}
\label{IALaplaceequation} v(r,z)=E_{r,z}\bigl[ A[r,T] \bigr]-
E_{r,z} \biggl[\int_r^T\psi
\bigl(s,Z_s,v(s,Z_s)\bigr) K(ds) \biggr]
\end{equation}
for $0\le r\le T$. According to Dynkin [(\citeyear{DynkinBranchingBook}),
Theorem 3.4.2], solutions to
(\ref{IALaplaceequation}) are unique when $E_{r,z}[ A[r,T] ]$ is
uniformly bounded in $z$ and $r\le T$, and hence in particular when
$A\in\cA_{(1)}$. Here we have the following result,
which will be proved in Section \ref{maxprinzSection}.

\begin{proposition}\label{UniquenessProp}For $A\in\cA^1_T$ and
$\psi$ as in (\ref{psiclass}), the function $v(r,z):=-\log\bE
_{r,\delta_z}[ e^{-J_A} ]$ is the unique finite and nonnegative
solution of (\ref{IALaplaceequation}).
\end{proposition}

\subsection{Statement of main results}\label{Statementsection}

Let $T>0$ be a fixed finite time horizon, $z\in{S}$ a fixed starting
point for $Z$ and $x_0\in\bR$ a given initial value. An
\textit{admissible strategy} will be a stochastic process $(x(t))_{0\le t\le
T}$ that is of the form $x(t)=x_0+\int_0^t\dot x(s)\,ds$ for an
integrable and $(\cF[0,t])$-progressively measurable process $(\dot
x(t))_{0\le t\le T}$. We also assume that $x(\cdot)$ satisfies the
fuel constraint
%
\begin{equation}
\label{liquidationconstraint}
x(T)=0,\qquad \mbox{$P_{0,z}$-a.s.}
\end{equation}
The \textit{cost} of an admissible strategy will be
%
\begin{equation}
\label{generalcostfct} E_{0,z} \biggl[ \int_0^T\bigl|
\dot x(s)\bigr|^p\eta(Z_s)\,ds+\int_{[0,T]}\bigl|x(s)\bigr|^p
A(ds) \biggr].
\end{equation}
Here, $p\in[2,\infty)$ and $\eta\dvtx {S}\to(0,\infty)$ is a measurable
function\footnote{We can always assume that $Z$ is of the form
$Z_t=(t,\widetilde Z_t)$, and so there is no loss of generality in assuming
the form $\eta(Z_t)$ rather than the form
$\eta(t,Z_t)$.} that will be further specified below. Our goal is the
minimization of this cost functional over all admissible strategies.

We are now ready to state our first main result, pertaining to the case
$\eta\equiv1$. It is actually a corollary of the more general
Theorem \ref{mainthm}, but the latter needs additional assumptions
and preparation. So we state this result here for the impatient reader.

\begin{theorem}\label{homogeneousthm}For $p\in[2,\infty)$ let $q$
be such that $\frac1p+\frac1q=1$ and take $A\in\cA^q_T$. Let $J_A$
be the corresponding $J$-functional of the superprocess with parameters
$Z$, $K(ds)=\frac1\beta\,ds$ and $\psi(t,z,\xi)=\xi^{1+\beta}$
for $\beta:=\frac1{p-1}$, and define the function
\[
v_\infty(r,z)=-\log\bE_{r,\delta_z}\bigl[ e^{-J_A} \Ind{
\{X_T=0\}} \bigr].
\]
Then
\[
x^*(t):=x_0\exp \biggl(-\int_0^tv_\infty(s,Z_s)^\beta
\,ds \biggr)
\]
is the unique admissible strategy that minimizes the cost functional
(\ref{generalcostfct}) for the choice $\eta\equiv1$. Moreover, the
minimal costs are given by
\[
E_{0,z} \biggl[ \int_0^T\bigl|\dot
x^*(s)\bigr|^p\,ds+\int_{[0,T]}\bigl|x^*(s)\bigr|^p
A(ds) \biggr]= |x_0|^pv_\infty(0,z).
\]
\end{theorem}

\begin{remark}In the proof of Theorem \ref{mainthm} it will be shown
that the optimal strategy is always below the linear strategy, that is,
%
\begin{equation}
\label{x*bound} \bigl|x^*(t)\bigr|\le|x_0| \frac{T-t}{T} \qquad\mbox{for $0\le t
\le T$; }
\end{equation}
see Remark \ref{x*boundRemark}.
\end{remark}

We now turn to extending the results from Theorem \ref{homogeneousthm}
to the minimization of the cost functional (\ref{generalcostfct}) with
nonconstant $\eta$.
Unless $\eta$ is bounded away from zero, this problem cannot be solved
by the standard class of superprocesses considered in Dynkin
(\citeyear{DynkinHistorical,DynkinBranchingBook}); we need the extended class of
superprocesses constructed in \citet{SchiedDiplompaper} by means of an
``$h$-transform.'' An analytical version of this transform was
found independently by \citet{EnglaenderPinsky}. Here we will use the
probabilistic version.
To make it work, we need the following assumption, which we will impose
from now on: for given $T>0$ there exists a constant $c_T>0$ such that
%
\begin{equation}
\label{almostharmonic}\quad \frac1{c_T}\eta(z)\le E_{r,z}\bigl[
\eta(Z_t) \bigr]\le c_T\eta(z) \qquad\mbox{for $0\le r\le t\le
T$ and $z\in{S}$.}
\end{equation}
We assume moreover that
%
\begin{equation}
\label{uniformalmostharmonic} E_{t,z}\bigl[ \eta(Z_T)
\bigr]\to\eta(z) \qquad\mbox{uniformly in $z$ as $t\ua T$.}
\end{equation}
By $\cM^\eta$ we denote the class of all nonnegative measures $\mu$ on
$({S},\cB)$ for which $\int\eta\,d\mu<\infty$ and by $B^+_\eta $ the
class of all bounded $\cB$-measurable functions $f\dvtx {S}\to\bR _+$ for
which there is a constant $c$ such that $f\le c\eta$. The proof of
the following result will be based on Schied
[(\citeyear{SchiedDiplompaper}), Theorem 2] and given in Section
\ref{h-trafosection}.

\begin{proposition}\label{h-transformexistProp}For $\beta\in(0,1]$
and under assumption (\ref{almostharmonic}), there exists an $\cM
^\eta$-valued Markov process $X=((X_t)_{t\le T},\cG(I),\bP_{r,\mu
})$, the superprocess with parameters $Z$, $K(ds)=\frac1{\beta}\eta
(Z_s) \,ds$ and $\psi(\xi)= (\frac\xi{\eta(z)} )^{1+\beta
}$, whose Laplace functionals are given by
%
\begin{equation}
\label{h-transformLaplacefunct} \bE_{r,\mu}\bigl[ e^{-\langle f,X_t\rangle}
\bigr]=e^{-\langle u(r,\cdot
),\mu\rangle},\qquad f\in B_\eta^+, \mu\in\cM^\eta, t\le
T,
\end{equation}
where $u$ is the unique solution in $ B_\eta^+$ of the integral equation
%
\begin{equation}
\label{h-transformLaplacefuncteq} u(r,z)=E_{r,z} \biggl[
f(Z_t)-\int_r^tu(s,Z_s)^{1+\beta}
\frac1{\beta \eta(Z_s)^\beta}\,ds \biggr].
\end{equation}
Moreover, to each $A\in\cA_T^1$ there exists a corresponding
$J$-functional, $J_A$, satisfying
%
\begin{equation}
\bE_{r,\mu}\bigl[ e^{-J_A} \bigr]=e^{-\langle v(r,\cdot),\mu\rangle},
\end{equation}
where $v$ solves
%
\begin{equation}
\label{hJALaplaceequation} v(r,z)=E_{r,z} \biggl[ A[r,T]-\int
_r^tv(s,Z_s)^{1+\beta}\frac1{
\beta \eta(Z_s)^\beta} \,ds \biggr].
\end{equation}
Furthermore, $v$ is the unique finite and nonnegative solution of
(\ref{hJALaplaceequation}).
\end{proposition}

\begin{example}[(Hyperbolic superprocesses)]
Let $Z$ be a one-dimensional Brownian motion stopped when first hitting
zero. Then one can take $\eta(z)=|z|^{1/\beta+\lambda}$ for some
$\lambda\in[0,1]$; see \citet{SchiedDiplompaper}, Example 1(i).
For $\beta=\lambda=1$, the corresponding superprocess was constructed
in \citet{FleischmannMueller}.
\end{example}

We are now ready to state our main result.

\begin{theorem}\label{mainthm}Suppose that $\eta$ is as above. For
$p\in[2,\infty)$ let $q$ be such that $\frac1p+\frac1q=1$, and take
$A\in\cA^1_T$ such that
%
\begin{equation}
\label{etaAassumption} \int_0^TE_{0,z}
\bigl[ \eta(Z_t)^{1-q}A[t,T]^q \bigr]\,dt<
\infty.
\end{equation}
For $\beta:=\frac1{p-1}$ let moreover $X=(X_t,\cG(I),\bP_{r,\mu})$
be the superprocess constructed in Proposition \ref
{h-transformexistProp}, and define
%
\begin{equation}
\label{vinftydef} v_\infty(r,y):=-\log\bE_{r,\delta_y}\bigl[
e^{-J_A} \Ind{\{X_T=0\}} \bigr]
\end{equation}
and
%
\begin{equation}
\label{x*} x^*(t):=x_0\exp \biggl(-\int_0^t
\frac{v_\infty(s,Z_s)^\beta}{\eta
(Z_s)^\beta}\,ds \biggr).
\end{equation}
Then $x^*$ is an admissible strategy and the unique minimizer of the
cost functional (\ref{generalcostfct}). Moreover, the minimal costs
are given by
\[
E_{0,z} \biggl[ \int_0^T\bigl|\dot
x^*(s)\bigr|^p\eta(Z_s)\,ds+\int_{[0,T]}\bigl|x^*(s)\bigr|^p
A(ds) \biggr]= |x_0|^pv_\infty(0,z).
\]
\end{theorem}

Note that the function $v_\infty(t,y)$ blows up as $t\ua T$. This fact
will create considerable difficulties. In fact, the most difficult part
in the proof of Theorem \ref{mainthm} will be to show that the
strategy $x^*$ defined through (\ref{x*}) is an admissible strategy
with finite cost. To prove this, we need sharp upper and lower bound
for the behavior of $v_\infty(t,y)$ as $t\ua T$. These estimates are
of independent interest and will be developed in Section \ref{auxSection}.

The blow-up of the function $v_\infty$ is of course linked to the fuel
constraint (\ref{liquidationconstraint}) required from admissible strategies.
A common question one therefore encounters
in relation to finite-fuel control problems is whether it may not be
reasonable to replace the sharp fuel constraint by a suitable
penalization term and thus avoid the singularity of the value
function.\vadjust{\goodbreak}
It turns out that in our context such a penalization approach can be
carried out without much additional effort.
To describe it, we define a \textit{relaxed strategy} as a stochastic
process $(x(t))_{0\le t\le T}$ that is of the form $x(t)=x_0+\int_0^t\dot x(s)\,ds$ for an integrable and $(\cF[0,t])$-progressively
measurable process $(\dot x(t))_{0\le t\le T}$.
The \textit{cost} of a relaxed strategy $x(\cdot)$ will be
%
\begin{equation}
\label{relaxedcostfct} \quad E_{0,z} \biggl[ \int_0^T\bigl|
\dot x(s)\bigr|^p\eta(Z_s)\,ds+\int_{[0,T]}\bigl|x(s)\bigr|^p
A(ds)+\varrho(Z_T)\bigl|x(T)\bigr|^p \biggr],
\end{equation}
where $p$, $A$ and $\eta$ are as above, and $\varrho\dvtx {S}\to\bR_+$
is a measurable penalty function that satisfies the following
assumptions for given $T$:
%
\begin{eqnarray}
\label{varrhocond} \varrho(y)&\le& c_\varrho\eta(y) \qquad\mbox{for a constant
$c_\varrho$ and all $y$;}
\nonumber\\[-8pt]\\[-8pt]
E_{t,Z_t}\bigl[ \varrho(Z_T) \bigr]&\lra&\varrho(Z_T)
\qquad\mbox{$P_{0,z}$-a.s. as $t\ua T$.}
\nonumber
\end{eqnarray}
By martingale convergence, the second condition is satisfied as soon as
$Z_T$ is $P_{0,z}$-a.s. measurable with respect to $\sigma(\bigcup_{t<T}\cF[0,t])$, which in turn holds when $Z$ is a Hunt process.
Moreover, it follows from (\ref{uniformalmostharmonic}) that both
conditions in (\ref{varrhocond}) are satisfied when $\varrho=c\eta$
for some constant $c\ge0$.

\begin{theorem}\label{relaxedthm}Suppose that $\eta$, $p$, $\beta$
and $X$ are as in Theorem \ref{mainthm} and that the measurable
penalty function $\varrho\dvtx {S}\to\bR_+$ satisfies (\ref{varrhocond}).
Let $A\in\cA_T^1$ be such that $A\{T\}=0$ $P_{0,z}$-a.s. For
%
\begin{equation}
v_\varrho(r,y)=-\log\bE_{r,\delta_y}\bigl[ e^{-J_A-\langle\varrho,X_T\rangle}
\bigr],
\end{equation}
the relaxed strategy
%
\begin{equation}
\label{xvarrho} x_\varrho(t):=x_0\exp \biggl(-\int
_0^t\frac{v_\varrho(s,Z_s)^\beta
}{\eta(Z_s)^\beta}\,ds \biggr)
\end{equation}
is the unique minimizer of the cost functional (\ref{relaxedcostfct})
in the class of relaxed strategies. Moreover, the minimal costs are
\begin{eqnarray*}
&&
E_{0,z} \biggl[ \int_0^T\bigl|\dot
x_\varrho(s)\bigr|^p\eta(Z_s)\,ds+\int
_{[0,T]}\bigl|x_\varrho(s)\bigr|^p A(ds)+
\varrho(Z_T)\bigl|x_\varrho(T)\bigr|^p \biggr]\\
&&\qquad=|x_0|^pv_\varrho(0,z).
\end{eqnarray*}
\end{theorem}

\section{Auxiliary results on superprocesses and their Laplace
functionals}\label{auxSection}

\subsection{A probabilistic version of the parabolic maximum
principle}\label{maxprinzSection}

Our first result is the
following proposition, which can be regarded as a probabilistic version
of a parabolic maximum principle for equations of the form (\ref
{IALaplaceequation}).

\begin{proposition}\label{MaxprinzProp}
Suppose that $A,\widetilde A\in\cA^1_T$. Let furthermore $ L$ and $
\widetilde L$
be nonnegative continuous additive functionals of $Z$. Suppose that
$A[r,t]\le\widetilde A[r,t]$ and $ L[r,t]\ge\widetilde L[r,t]$ for
$0\le r\le t\le
T$. For some $q\ge1$ let $v$ and $\widetilde v$ be finite and nonnegative
solutions of the integral equations
\begin{eqnarray*}
v(r,{z})&=&E_{r,{z}} \biggl[ A[r,T]-\int_r^Tv(t,{Z}_t)^q
L(dt) \biggr],
\\
\widetilde v(r,{z})&=&E_{r,{z}} \biggl[ \widetilde A[r,T]-\int
_r^T\widetilde v(t,{Z}_t)^q
\widetilde L(dt) \biggr].
\end{eqnarray*}
Then we have $v\le\widetilde v$.
\end{proposition}

For this and other proofs we will need the following special version of
the general Feynman--Kac formula from \citet{DynkinBranchingBook},
Theorem 4.1.1.

\begin{proposition}[{[Dynkin (\citeyear{DynkinBranchingBook})]}]\label{Feynman-Kac}Suppose that $B$
is a signed additive functional of $Z$ whose total variation belongs to
$\cA^1_T$. Let furthermore $C$ be a continuous and nonnegative
additive functional of $Z$. Let
\[
g(r,z):=E_{r,z} \biggl[ \int_{[r,T]}e^{-C[r,s]}
B(ds) \biggr].
\]
When $E_{r,z}[ \int_r^tg(s,Z_s) C(ds) ]$ is well defined and finite
for all $(r,z)$, then $g$ is the unique solution of the linear integral equation
%
\begin{equation}
g(r,z)=E_{r,z} \bigl[ B[r,T] \bigr]-E_{r,z} \biggl[ \int
_r^tg(s,Z_s) C(ds) \biggr].
\end{equation}
\end{proposition}

\begin{pf*}{Proof of Proposition \ref{MaxprinzProp}}
Via the outer regularity of finite Borel measures on $[0,T]$, our
condition $\widetilde L[r,t]\le L[r,t]$ implies that $ \widetilde
L(dt)=\varphi _t L(dt)$ for some $[0,1]$-valued function $\varphi_t$,
which can be chosen to be progressively measurable. One sees that
$u:=\widetilde v-v$ satisfies
\begin{eqnarray*}
u(r,{z})&=&E_{r,{z}} \biggl[ \widetilde A[r,T]-A[r,T]+\int
_r^Tv(t,Z_t)^q(1-
\varphi_t) L(dt)\\
&&\hspace*{108pt}{}-\int_r^Tu(t,{Z}_t)w(t,Z_t)
\widetilde L(dt) \biggr],
\end{eqnarray*}
where
\[
w(t,z):=\cases{\displaystyle\frac{\widetilde
v(t,z)^q-v(t,z)^q}{\widetilde
v(t,z)-v(t,z)}, &\quad if $\widetilde v(t,z)\neq
v(t,z)$,
\vspace*{2pt}\cr
0, &\quad otherwise.}
\]
Note that $w(t,z)$ is nonnegative. We define an additive functional $B$
of $Z$ via
\[
B[r,t]:=\widetilde A[r,t]-A[r,t]+\int_r^tv(s,Z_s)^q(1-
\varphi_s) L(ds).
\]
It is nonnegative and belongs to $\cA^1_T$ because
\[
E_{r,z} \biggl[ \int_r^Tv(s,Z_s)^q(1-
\varphi_s) L(ds) \biggr]\le E_{r,z} \biggl[ \int
_r^Tv(s,Z_s)^q L(ds)
\biggr],
\]
and the expectation on the right is finite by assumption. Moreover,
$C(dt):=w(t,Z_t) \widetilde L(dt)$ is a continuous and nonnegative additive
functional. We have
\begin{eqnarray*}
&&
E_{r,z} \biggl[ \int_r^T\bigl|u(s,Z_s)\bigr|
C(ds) \biggr]
\\
&&\qquad\leq E_{r,{z}} \biggl[ \int_r^T
\widetilde v(t,{Z}_t)^q\Ind{\bigl\{ w(t,Z_t)
\neq0\bigr\}} \widetilde L(dt) \biggr]\\
&&\qquad\quad{}+E_{r,{z}} \biggl[ \int
_r^T v(t,{Z}_t)^q\Ind {
\bigl\{ w(t,Z_t)\neq0\bigr\}} \widetilde L(dt) \biggr].
\end{eqnarray*}
Due to our assumptions, both expectations on the right are finite, and
so $E_{r,z} [ \int_r^Tu(s,Z_s) C(ds) ]$ is well defined and finite
for all $(r,z)$. Therefore, Proposition \ref{Feynman-Kac} yields
that
\[
u(r,z)=E_{r,z} \biggl[ \int_{[r,T]}e^{-C[r,s]}
B(ds) \biggr],
\]
which is nonnegative. Hence, $\widetilde v\ge v$.
\end{pf*}

\begin{pf*}{Proof of Proposition \ref{UniquenessProp}}
By Dynkin [(\citeyear{DynkinHistorical}), Theorem 1.2], we have
$\bE_{r,\delta_z}[ J_A ]\le E_{r,z}[ A[r,T] ]$, which is finite for all
$(r,z)$ since $A\in \cA_T^1$. Hence $J_A<\infty$
$\bP_{r,\delta_z}$-a.s. and so $v(r,z)$ is finite for all $(r,z)$. Note
that
%
\begin{eqnarray}
\label{Expwelldefauxeq} && E_{r,z} \biggl[\int_r^Tv(s,Z_s)^{1+\beta}a(s,Z_s)
K(ds) \biggr]
\nonumber\\[-8pt]\\[-8pt]
&&\qquad =E_{r,z} \biggl[\int_r^T\psi
\bigl(s,Z_s,v(s,Z_s)\bigr) K(ds) \biggr]<\infty\nonumber
\end{eqnarray}
for all $(r,z)$. When $\widetilde v$ is another finite and nonnegative
solution of (\ref{IALaplaceequation}), then
(\ref{Expwelldefauxeq}) holds also for $\widetilde v$. Therefore
we may
apply Proposition \ref{MaxprinzProp} with $L(ds)=\widetilde
L(ds)=a(s,Z_s)
K(ds)$ and $q=1+\beta$ to get $v\le\widetilde v$. Interchanging the
roles of
$v$ and $\widetilde v$ yields the uniqueness of solutions.
\end{pf*}

\begin{example}[(Laplace functionals for the total mass process)]
\label{totalmassexample}
Consider the superprocess with parameters $Z$, $K(ds)=\gamma\,ds$ for
a constant $\gamma>0$, and $\psi(s,z,\xi)=\xi^{1+\beta}$ for
$\beta\in(0,1]$.
Let $\nu$ be a finite and nonnegative Borel measure on $[0,T]$.
Clearly, $\nu$ can be regarded as an element of $\cA^1_T$. The
corresponding $J$-functional is given by
%
\begin{equation}
\label{Inu} J_\nu=\int\langle1,X_t\rangle\nu(dt)
\end{equation}
as can easily be seen from (\ref{IAdef}). Its log-Laplace functional
\[
v(r,z):=-\log\bE_{r,\delta_z} \bigl[ e^{-J_\nu} \bigr]=-\log\bE
_{r,\delta_z} \bigl[ e^{-\int_{[r,T]}\langle1,X_t\rangle\nu
(dt)} \bigr]
\]
is in fact independent of $z$. Indeed, when $\nu=\lambda\delta_t$
for some $\lambda\ge0$ and $t\in[r,T]$, then $v(r,z)=0$ for $r>t$ and
%
\begin{equation}
\label{totalmasshomLaplaceeq} v(r,z)=-\log\bE \bigl[ e^{-\lambda\langle1,X_t\rangle} \bigr]=
\frac\lambda{\bigl(1+\gamma\beta(t-r)\lambda^{\beta}\bigr)^{1/\beta}}
\end{equation}
for $r\le t$ and $\lambda\ge0$, as can be shown by a straightforward
computation based on the integral equation (\ref
{vlog-Laplaceequation}). When $\nu$ is a positive linear combination of Dirac
measures, then we can use the Markov property of $X$ to conclude that
$v(r,z)$ is independent of $z$. For general $\nu$ we use an
approximation argument.
Alternatively, one can use the fact that, for superprocesses with
homogeneous branching, the total mass process, $\langle1,X_t\rangle$,
is itself a one-dimensional Markov process.
It follows that the function $v$ is the unique nonnegative solution of
the integral equation
\[
v(r)=\nu[r,T]-\gamma\int_r^Tv(s)^{1+\beta}
\,ds.
\]
\end{example}

\subsection{Estimates for the Laplace functionals of
$J$-functionals}\label{HomogeneousLaplaceestimatesSection}

Throughout this section, let $X=(X_t,\cG(I),\bP_{r,\mu})$ be the
superprocess with one-particle motion $Z$, $K(ds)=\gamma\,ds$ for a
constant $\gamma>0$ and $\psi(s,z,\xi)=\xi^{1+\beta}$ for $\beta
\in(0,1]$.
A key ingredient in the proof of Theorem \ref{mainthm} will be
inequality (\ref{compeq}) in the following theorem. This inequality
gives a bound on the Laplace transform of a $J$-functional and is also
of independent interest. It extends the upper bound of the following
estimate for the Laplace functionals of $X$ from Section 5 of
\citet{SchiedDisspaper}:
%
\begin{equation}
\label{Laplaceineq} E_{r,z}\bigl[ V_{T-r}f(Z_T)
\bigr]\le-\log\bE_{r,\delta_z}\bigl[ e^{-\langle
f,X_T\rangle} \bigr]\le
V_{T-r}E_{r,z}\bigl[ f(Z_T) \bigr]
\end{equation}
for $r<T$ and $f\ge0$, where $V_t\dvtx \bR_+\to\bR_+$ denotes the
nonlinear semigroup
\[
V_ty=\frac y{(1+\gamma\beta t y^\beta
)^{1/\beta}},\qquad y\ge0, t\ge0.
\]
For $\beta=1$, estimate (\ref{Laplaceineq}) can actually be extended
to functions $f$ with arbitrary sign. The lower bound in (\ref
{Laplaceineq}) will be needed in the proof of Theorem \ref{relaxedthm}.
Note that, by means of Example \ref{totalmassexample}, the following
inequality (\ref{compeq}) coincides with the upper bound in (\ref
{Laplaceineq}) for the additive functional $A(dt)=f(Z_t) \delta_T(dt)$.

\begin{theorem}\label{condexpthm}For the superprocess with
homogeneous branching rate $\gamma$, let $J_A$ be the $J$-functional
associated with a given $A\in\cA^1_T$. For $r\le T$ and $z\in{S}$
fixed, define moreover the finite and nonnegative Borel measure $\alpha
_{r,z}(ds)$ on $[r,T]$ by
\[
\int f(s) \alpha_{r,z}(ds)=E_{r,z} \biggl[ \int
_{[r,T]}f(s) A(ds) \biggr]
\]
for bounded measurable $f\dvtx [r,T]\to\bR$. Then the conditional
expectation of $J_A$ given the evolution of the total mass process,
$(\langle1,X_t\rangle)_{r\le t\le T}$, is given by
%
\begin{equation}
\label{condexpeq} \bE_{r,\delta_z} \bigl[ J_A |
\langle1,X_t\rangle, r\le t\le T \bigr]=\int_{[r,T]}
\langle1,X_t\rangle\alpha_{r,z}(dt).
\end{equation}
Moreover,
%
\begin{equation}
\label{compeq} \bE_{r,\delta_z}\bigl[ e^{-J_A} \bigr]\ge
\bE_{r,\delta_z} \bigl[ e^{-\int
_{[r,T]}\langle1,X_t\rangle\alpha_{r,z}(dt)} \bigr].
\end{equation}
\end{theorem}
\begin{pf} To prove (\ref{condexpeq}), we can assume without loss
of generality that $A$ is bounded. For $\lambda\ge0$ and a bounded
nonnegative Borel measure $\mu$ on $[0,T]$, let
\[
v_\lambda(s,y):=-\log\bE_{s,\delta_y}\bigl[ e^{-\lambda J_A-\int
\langle1,X_t\rangle\mu(dt)} \bigr].
\]
Then $v_\lambda$ is the unique nonnegative solution of
%
\begin{equation}
\label{vsequation}\quad v_\lambda(s,y)=E_{s,y}\bigl[ \lambda A[s,T]
\bigr]+\mu[s,T]-\gamma\int_s^TE_{s,y}
\bigl[ v_\lambda(t,Z_t)^{1+\beta} \bigr]\,dt.
\end{equation}
It follows that
%
\begin{equation}
\label{wdefandbddeq}\qquad \bE_{r,\delta_z}\bigl[ J_A e^{-\int\langle1,X_t\rangle\mu(dt)}
\bigr]=-\frac{d}{d\lambda} \bigg|_{\lambda=0}e^{-v_\lambda
(r,z)}=e^{-v_0(r,z)}\,
\frac{\partial v_\lambda(r,z)}{\partial
\lambda} \bigg|_{\lambda=0}.
\end{equation}
For $\lambda=0$, we have
\[
v_0(s,y)=-\log\bE_{s,\delta_y} \bigl[ e^{-\int\langle1,X_t\rangle
\mu(dt)} \bigr],
\]
which is in fact independent of $z$ due to the assumed homogeneity of
the branching mechanism; see Example \ref{totalmassexample}.
Hence, $w:=\partial v_\lambda/\partial\lambda|_{\lambda=0}$ solves
\[
w(s,y)=E_{s,y}\bigl[ A[s,T] \bigr]-\gamma(1+\beta)\int
_s^TE_{s,y}\bigl[
w(t,Z_t) \bigr]v_0(t)^{\beta} \,dt.
\]
Here, interchanging differentiation and integration is justified due to
the uniform boundedness of $A$ and, hence, of $v_\lambda$ and $w$ [the
boundedness of the latter being implied by (\ref{wdefandbddeq})].
Due to the general Feynman--Kac formula of Proposition
\ref{Feynman-Kac}, $w$~is given by
\begin{eqnarray*}
w(r,z)&=&E_{r,z} \biggl[ \int_{[r,T]}e^{-\gamma(1+\beta)\int
_r^tv_0(u)^\beta \,du}
A(dt) \biggr]
\\
&=&\int_{[r,T]}e^{-\gamma(1+\beta)\int_r^tv_0(u)^\beta \,du} \alpha _{r,z}(dt).
\end{eqnarray*}
Therefore
%
\begin{equation}
\label{condexp1}\quad \bE_{r,\delta_z}\bigl[ J_A e^{-\int\langle1,X_t\rangle\mu(dt)}
\bigr]=e^{-v_0(r)}\int_{[r,T]}e^{-\gamma(1+\beta)\int_r^tv_0(u)^\beta
\,du}
\alpha_{r,z}(dt).
\end{equation}
Now we take $A(dt)=\nu(dt)$ for a nonnegative finite Borel measure
$\nu$
on $[0,T]$ and recall from (\ref{Inu}) that in this case $J_A=J_\nu
=\int\langle1,X_t\rangle\nu(dt)$. We then get
%
\begin{equation}
\label{condexp2} \bE_{r,\delta_z} \bigl[ J_\nu e^{-\int\langle1,X_t\rangle\mu
(dt)}
\bigr]=e^{-v_0(r)}\int_{[r,T]}e^{-\gamma(1+\beta)\int
_r^tv_0(u)^\beta \,du} \nu(dt).
\end{equation}
Comparing (\ref{condexp1}) with (\ref{condexp2}) and recalling
(\ref{Inu}) yields
\begin{eqnarray*}
\bE_{r,\delta_z}\bigl[ J_A e^{-\int\langle1,X_t\rangle\mu(dt)}
\bigr]&=&
\bE_{r,\delta_z} \bigl[ J_{\alpha_{r,z}} e^{-\int\langle
1,X_t\rangle\mu(dt)} \bigr]
\\
&=&\bE_{r,\delta_z} \biggl[ \int_{[r,T]}
\langle1,X_t\rangle\alpha _{r,z}(dt)e^{-\int\langle1,X_t\rangle\mu(dt)}
\biggr].
\end{eqnarray*}
Varying $\mu$ and applying a monotone class argument yields (\ref
{condexpeq}).

Now we prove (\ref{compeq}). We have
\begin{eqnarray*}
\bE_{r,\delta_z}\bigl[ e^{-J_A} \bigr]&=&\bE_{r,\delta
_z} \bigl[
\bE_{r,\delta_z}\bigl[ e^{-J_A} | \langle1,X_t\rangle, r\le
t\le T \bigr] \bigr]
\\
&\ge&\bE_{r,\delta_z}\bigl[ e^{-\bE_{r,\delta_z}[ J_A | \langle
1,X_t\rangle, r\le t\le T ]} \bigr]
\\
&=&\bE_{r,\delta_z} \bigl[ e^{-\int_{[r,T]}\langle1,X_t\rangle
\alpha_{r,z}(dt)} \bigr],
\end{eqnarray*}
where we have used Jensen's inequality for conditional expectations in
the second step and (\ref{condexpeq}) in the third.
\end{pf}

Recall that $K(ds)=\gamma \,ds$ and $\psi(s,z,\xi)=\xi^{1+\beta}$.
Let us also mention that the following estimates will be extended to
the case of nonhomogeneous branching in Proposition
\ref{inhomLaplaceestimatesProp}.

\begin{proposition}\label{kLaplaceestimatesProp}Let $J_A$ be the
$J$-functional associated with $A\in\cA^1_T$. Then, for $k\ge0$ and $r<T$,
%
\begin{equation}
\label{log-Laplaceestimatewithk}\qquad -\log\bE_{r,\delta_z}\bigl[ e^{-J_A-k\langle1,X_T\rangle}
\bigr]\le E_{r,z}\bigl[ A[r,T] \bigr]+\frac k{(1+\gamma
\beta(T-r)k^{\beta})^{1/\beta}}
\end{equation}
and
%
\begin{equation}
\label{log-Laplaceestimatewithinfinity} -\log\bE_{r,\delta_z}\bigl[
e^{-J_A} \Ind{\{X_T=0\}} \bigr]\le E_{r,z}\bigl[
A[r,T] \bigr]+\frac1{(\gamma\beta(T-r))^{1/\beta}}.
\end{equation}
\end{proposition}
\begin{pf} We can assume without loss of generality that $A$ is
bounded. From (\ref{compeq}) we get
%
\begin{equation}
\label{logLaplaceauxestimate}\qquad-\log\bE_{r,\delta
_z}\bigl[ e^{-J_A-k\langle1,X_T\rangle}
\bigr]\le-\log\bE_{r,\delta
_z} \bigl[ e^{-\int_{[r,T]}\langle1,X_t\rangle\mu(dt)-k\langle
1,X_T\rangle} \bigr],
\end{equation}
where $\mu([s,t])=E_{r,z}[ A[s,t] ]$ for $r\le s\le t\le T$. As
noted in Example \ref{totalmassexample}, the right-hand side of
(\ref{logLaplaceauxestimate}) is independent of $z$ and equal to
$v(r)$, where $v$ solves the integral equation
\[
v(t)=k+\mu\bigl([t,T]\bigr)-\gamma\int_t^Tv(s)^{1+\beta}
\,ds,\qquad 0\le t\le T.
\]
Assertion (\ref{log-Laplaceestimatewithk}) now follows from an
application of Lemma \ref{Gronwall-stylelemma}, which is stated
below. Inequality (\ref{log-Laplaceestimatewithinfinity}) is
obtained by sending $k$ to infinity in (\ref{log-Laplaceestimatewithk}).
\end{pf}

\begin{lemma}\label{Gronwall-stylelemma}Suppose that $a\dvtx [0,T]\to\bR
_+$ is a measurable function, $k\ge0$ is a constant and $v\dvtx [0,T]\to
\bR_+$ solves the integral equation
\[
v(r)=k+a(r)-\gamma\int_r^Tv(s)^{1+\beta}
\,ds,\qquad 0\le r\le T.
\]
Then
%
\begin{equation}
v(t)\le a(t)+\frac k{(1+\gamma\beta(T-t)k^{\beta}
)^{1/\beta
}},\qquad 0\le t\le T.
\end{equation}
\end{lemma}
\begin{pf} The function
\[
u(t):=\frac k{(1+\gamma\beta(T-t)k^{\beta})^{1/\beta}}
\]
satisfies $u(T)=k$ and solves
\[
u(r)=k-\gamma\int_r^Tu(s)^{1+\beta} \,ds.
\]
Let $\widetilde v(t):=v(T-t)$, and define $\widetilde u$ and
$\widetilde a$ accordingly. The
function $w(t):=\widetilde v(t)-\widetilde u(t)-\widetilde a(t)$ is
absolutely continuous
and satisfies for a.e. $t$
\[
w'(t)=-\gamma \bigl(\widetilde v(t)^{1+\beta
}-\widetilde
u(t)^{1+\beta} \bigr)=-\gamma \bigl(\widetilde v(t)-\widetilde u(t) \bigr)
f(t),
\]
where
\[
f(t)=\cases{\displaystyle  \frac{\widetilde v(t)^{1+\beta
}-\widetilde u(t)^{1+\beta
}}{\widetilde v(t)-\widetilde u(t)}, &\quad for $\widetilde v(t)\neq\widetilde u(t) $,
\vspace*{2pt}\cr
0, &\quad otherwise.}
\]
Since $f\ge0$ and $a\ge0$, it follows that $w'(t)\le-\gamma
w(t)f(t)$ for a.e. $t\in[0,T]$. When letting $w_0(t):=e^{-\gamma\int
_0^tf(s) \,ds}$, we have
\[
\biggl(\frac{w(t)}{w_0(t)} \biggr)'=\frac
{w'(t)w_0(t)-w(t)w_0'(t)}{w_0(t)^2}\le0
\]
and so
\[
\frac{w(t)}{w_0(t)}\le\frac{w(0)}{w_0(0)}=v(T)-u(T)-a(T)=0.
\]
It follows that $w(t)\le0$ and in turn that $v\le a+u$.
\end{pf}

\subsection{An ``$h$-transform'' for
superprocesses}\label{h-trafosection}

In this section, we prove Proposition \ref{h-transformexistProp} and
extend the estimates from Proposition \ref{kLaplaceestimatesProp}
to certain superprocesses with inhomogeneous branching characteristics.
Our approach is based on the ``$h$-transform'' for
superprocesses that was introduced independently by \citet
{EnglaenderPinsky} and \citet{SchiedDiplompaper}. Whereas the first
approach is primarily analytical, the latter approach is probabilistic,
and it is the one we are going to use here. It is based on the
following space--time harmonic function of $Z$:
%
\begin{equation}
\label{hdef} h(r,z):=E_{r,z}\bigl[ \eta(Z_T) \bigr],\qquad
\mbox{$0\le r\le T$, $z\in{S}$.}
\end{equation}
We define the function $\psi(z,\xi)= (\frac\xi{\eta(z)}
)^{1+\beta}$, $\xi\ge0$, $z\in{S}$ and a continuous nonnegative
additive functional $K$ of $Z$ by $K(dt)=\frac1\beta\eta(Z_t) \,dt$.
Then we have $\widetilde\psi(z,\xi):=\psi(z,\eta(z)\xi)=\xi
^{1+\beta
}$. Moreover, by (\ref{almostharmonic}),
\[
E_{r,z} \bigl[ K[r,t] \bigr]=\frac1\beta\int_r^tE_{r,z}
\bigl[ \eta (Z_t) \bigr] \,dt\le\frac1\beta c_T(t-r)\eta(z)
\]
for $0\le r\le t\le T$ and $z\in{S}$.
Therefore, both $\psi$ and $K$ satisfy the conditions of Theorem 2
from \citet{SchiedDiplompaper}, which hence implies the existence of a
$(Z,K,\psi)$-superprocess $X$, defined here up to time $T$, for which
the function $u$ from (\ref{h-transformLaplacefunct}) uniquely solves
\begin{eqnarray*}
u(r,z)&=&E_{r,z} \biggl[ f(Z_t)-\int_r^t
\psi\bigl(Z_s,u(s,Z_s)\bigr) K(ds) \biggr]
\\
&=&E_{r,z} \biggl[ f(Z_t)-\int_r^tu(s,Z_s)^{1+\beta}
\frac1{\beta\eta (Z_s)^\beta} \,ds \biggr].
\end{eqnarray*}
This implies the first part in the assertion of Proposition
\ref{h-transformexistProp}.

To prove the remaining part of Proposition \ref{h-transformexistProp}
and to prepare for the proof of Theorem \ref{mainthm}, we need
to recall the construction of $X$ given in \citet{SchiedDiplompaper}.
One first introduces Doob's $h$-transform of the process $Z$, that is,
the Markov process $Z^h=(Z_t,\cF(I),P^h_{r,z})$\vadjust{\goodbreak} (defined up to the
time horizon $T$) where
\[
P^h_{r,z}[ A ]=\frac1{h(r,z)}E_{r,z}\bigl[
\eta(Z_T)\Ind{A} \bigr],\qquad A\in\cF[r,T]
\]
[note that $h>0$ by (\ref{almostharmonic})]. For any additive
functional $B$ of $Z$ one then defines an additive functional $B_h$ of
$Z^h$ by
\[
B_h(ds)=\frac1{h(s,Z_s)} B(ds).
\]
By the right property of $Z$, the process $h(s,Z_s)$ is a right
continuous $P_{r,z}$-martingale. Hence, $h(s,Z_s)$ is equal to the
optional projection of the constant process $t\mapsto\eta(Z_T)$.
Therefore, for $0\le r\le t\le T$,
%
\begin{eqnarray}
\label{BvsBh} E_{r,z} \bigl[ B[r,t] \bigr]&=&E_{r,z} \biggl[
\int_{[r,t]}h(s,Z_s) B_h(ds)
\biggr]=E_{r,z} \bigl[ B_h[r,t] \eta(Z_T)
\bigr]
\nonumber\\[-8pt]\\[-8pt]
&=&h(r,z) E^h_{r,z}\bigl[ B_h[r,t] \bigr],\nonumber
\end{eqnarray}
where we have used Theorem 57 in Chapter VI of \citet
{DellacherieMeyerB} in the second step.

With this notation,
\[
K_h(ds)=\frac1{h(s,Z_s)} K(ds)=\frac{\eta(Z_s)}{\beta h(s,Z_s)} \,ds
\]
is a bounded and continuous additive functional of $Z^h$, and the function
\[
\psi_h(t,z,\xi)=\psi\bigl(z,h(t,z)\xi\bigr)
\]
is of the form (\ref{psiclass}). Therefore, up to the time horizon
$T$, we can define the $(Z^h,K_h,\psi_h)$-superprocess $X^h=(X_t^h,\cG
(I),\bP_{r,\mu}^h)$, for example, via Theorem 1.1 in \citet
{DynkinHistorical}. The $(Z,K,\psi)$-superprocess $X$ under $\bP
_{r,\mu}$ is then defined as the law of
%
\begin{equation}
\label{Xhtrafo} X_t(dz):=\frac1{h(t,z)} X^h_t(dz)
\end{equation}
under\vspace*{1pt} $\bP_{r,h.\mu}^h$, where $h.\mu$ denotes the measure $h(r,z)
\mu(dz)$. Using once again Theorem 57 in Chapter VI of \citet
{DellacherieMeyerB}, one checks that the log-Laplace functionals of $X$
are indeed given by (\ref{h-transformLaplacefunct}), (\ref
{h-transformLaplacefuncteq}).

\begin{pf*}{Proof of the second part of Proposition
\ref{h-transformexistProp}} For $A\in\cA_{(1)}\cap\cA_T^1$ it is
clear from (\ref{IAdef}) and (\ref{Xhtrafo}) that $J_A$ must be
defined as the $J$-functional $J_{A_h}$ for~$X^h$, and this
identification carries over to all $A\in\cA_T^1$ by approximation.
As above, one then checks that
\begin{eqnarray*}
v(r,z):\!&=&-\log\bE_{r,\delta_z}\bigl[ e^{-J_A} \bigr]=-\log\bE
^h_{r,h(r,z)\delta_z}\bigl[ e^{-J_{A_h}} \bigr]
\\
&=&-h(r,z)\log\bE^h_{r,\delta_z}\bigl[ e^{-J_{A_h}} \bigr]
\end{eqnarray*}
solves (\ref{hJALaplaceequation}).\vadjust{\goodbreak}

Conversely, when $A\in\cA_T^1$ and $\widetilde v$ is a nonnegative solution
of (\ref{hJALaplaceequation}), then $\widetilde
v^h(r,z):=\widetilde{v}(r,z)/h(r,z)$ solves
%
\begin{equation}
\label{wtvheq}\qquad \widetilde v^h(r,z)=E_{r,z}^h
\bigl[ A_h[r,T] \bigr]-E_{r,z}^h \biggl[ \int
_r^T\psi _h\bigl(s,Z_s,
\widetilde v^h(s,Z_s)\bigr) K_h(ds)
\biggr].
\end{equation}
By\vspace*{1pt} (\ref{BvsBh}), $A_h$ belongs to the class $\cA_T^1$ for $Z^h$,
and so Proposition \ref{UniquenessProp} implies that $\widetilde
v^h$ is the
unique finite and nonnegative solution of (\ref{wtvheq}). But this
equation is also solved by $v^h(r,z)=h(r,z)v(r,z)$, which gives the
uniqueness of solutions to the equation (\ref{hJALaplaceequation}).
\end{pf*}

Now we turn toward generalizing the results from Section
\ref{HomogeneousLaplaceestimatesSection} to superprocesses with
inhomogeneous, state-dependent branching mechanism as constructed in
Proposition \ref{h-transformexistProp}.

\begin{proposition}\label{inhomLaplaceestimatesProp}Let $X$ be the
superprocess constructed in Proposition \ref{h-transformexistProp},
and let $J_A$ be the $J$-functional associated with $A\in\cA^1_T$. Then
%
\begin{equation}
\label{log-Laplaceestimatewithketa}\qquad -\log\bE_{r,\delta_z}\bigl[
e^{-J_A-k\langle\eta,X_T\rangle} \bigr]\le E_{r,z}\bigl[ A[r,T] \bigr]+
\frac{ h(r,z)k}{(1+c_T^{-\beta}(T-r)k^{\beta
})^{1/\beta}}
\end{equation}
and
%
\begin{equation}
\label{log-Laplaceestimatewithinfinityeta} -\log\bE_{r,\delta_z}\bigl[
e^{-J_A} \Ind{\{X_T=0\}} \bigr]\le E_{r,z}\bigl[
A[r,T] \bigr]+\frac{c_Th(r,z)}{(T-r)^{1/\beta}},
\end{equation}
where $c_T$ is the constant from (\ref{almostharmonic}).
\end{proposition}
\begin{pf}With the notation introduced above, we have
\[
-\log\bE_{r,\delta_z}\bigl[ e^{-J_A-k\langle\eta,X_T\rangle} \bigr]=-\log
\bE^h_{r,h(r,z)\delta_z}\bigl[ e^{-J_{A_h}-k\langle1,X^h_T\rangle} \bigr]=h(r,z)v^h(r,z),
\]
where $v^h(r,z)$ solves
%
\begin{eqnarray}
\label{vhequationA} v^h(r,z)&=&k+E^h_{r,z}
\bigl[ A_h[r,T] \bigr]-E^h_{r,z} \biggl[ \int
_r^T\psi_h\bigl(t,Z_t,v^h(t,Z_t)
\bigr) K_h(dt) \biggr]\hspace*{-30pt}
\nonumber\\[-8pt]\\[-8pt]
&=&k+E^h_{r,z} \bigl[ A_h[r,T]
\bigr]-E^h_{r,z} \biggl[ \int_r^T
\frac 1\beta \biggl(\frac{h(t,Z_t)}{\eta(Z_t)} \biggr)^\beta v^h(t,Z_t)^{1+\beta}
\,dt \biggr].\hspace*{-30pt}\nonumber
\end{eqnarray}
By (\ref{almostharmonic}), $\frac{h(t,Z_t)}{\eta(Z_t)}\ge\frac
1{c_T}$. Applying Proposition \ref{MaxprinzProp} with $L(dt)=\frac
1\beta(\frac{h(t,Z_t)}{\eta(Z_t)})^\beta$ and $\widetilde
L(dt):=\beta
^{-1}c_T^{-\beta} \,dt$ yields that
\[
v^h(r,z)\le\widetilde v^h(r,z):=-\log\widetilde
\bE_{r,\delta
_z}\bigl[ e^{-J_{A_h}-k\langle1,\widetilde X_T\rangle} \bigr],
\]
where $\widetilde X=(\widetilde X_t,\cG(I),\widetilde\bP_{r,\mu
})$ is the superprocess
with one-particle motion $Z^h$, branching function $\widetilde\psi
(\xi)=\xi
^{1+\beta}$ and branching\vadjust{\goodbreak} functional $\widetilde L$. Proposition
\ref{kLaplaceestimatesProp} yields that
\begin{eqnarray*}
\widetilde v^h(r,z)&\le& E^h_{r,z} \bigl[
A_h[r,T] \bigr]+\frac
{k}{(1+\beta
^{-1}c_T^{-\beta}\beta(T-r)k^{\beta})^{1/\beta}}
\\
&=&\frac1{h(r,z)}E_{r,z} \bigl[ A[r,T] \bigr]+\frac
{k}{(1+c_T^{-\beta} (T-r)k^{\beta})^{1/\beta}}.
\end{eqnarray*}
This proves (\ref{log-Laplaceestimatewithketa}). Sending $k$ to
infinity gives (\ref{log-Laplaceestimatewithinfinityeta}).
\end{pf}

We also need a lower bound in case $A=0$. To this end, we define
%
\begin{equation}
\label{crTdefinition} c_{r,T}:=\sup_{t\in[r,T]}\sup
_z\frac{h(t,z)}{\eta(z)}.
\end{equation}
It follows from (\ref{almostharmonic}) that $c_{r,T}$ is finite for
all $r\in[0,T]$ and from (\ref{uniformalmostharmonic}) that
$c_{r,T}\searrow1$ as $r\ua T$.

\begin{lemma}\label{lowerLaplaceestLemma}$\!\!\!$Let $X$ be the
superprocess constructed in Proposition \ref{h-transformexistProp}.~Then
\[
-\log\bP_{r,\delta_z}[ X_T=0 ]\ge\frac
{h(r,z)}{c_{r,T}(T-r)^{1/\beta}}.
\]
\end{lemma}
\begin{pf}We have
\begin{eqnarray*}
-\log\bP_{r,\delta_z}[ X_T=0 ]&=&-\lim_{k\ua
\infty
}
\log\bE_{r,\delta_z}\bigl[ e^{-k\langle\eta,X_T\rangle} \bigr]
\\
&=&-\lim_{k\ua\infty} \log\bE_{r,h(r,z)\delta_z}^h\bigl[
e^{-k\langle
1,X^h_T\rangle} \bigr]
\\
&=&h(r,z)\lim_{k\ua\infty}v_k^h(r,z),
\end{eqnarray*}
where $v^h_k$ solves (\ref{vhequationA}) for $A_h=0$. Using $\frac
{h(t,Z_t)}{\eta(Z_t)}\le{c_{r,T}}$ for $r\le t\le T$ and applying
Proposition \ref{MaxprinzProp} with $ L(dt)=\beta^{-1}c_{r,T}^\beta
\,dt$ and $\widetilde L(dt)=\frac1\beta(\frac{h(t,Z_t)}{\eta
(Z_t)})^\beta
$ hence yields that
\[
v_k^h(s,z)\ge-\log\widehat\bE_{s,\delta_z}\bigl[
e^{-k\langle
1,\widehat
X_T\rangle} \bigr],\qquad r\le s\le T,
\]
where $\widehat X=(\widehat X_t,\cG(I),\widehat\bP_{s,\mu})$ is
the superprocess
with one-particle motion $Z^h$, branching function $\widehat\psi(\xi
)=\xi
^{1+\beta}$ and branching functional $L$. By (\ref{totalmasshomLaplaceeq}),
\[
-\log\widehat\bE_{s,\delta_z}\bigl[ e^{-k\langle1,\widehat
X_T\rangle} \bigr]=\frac k{
(1+c_{r,T}^\beta(T-s)k^{\beta})^{1/\beta}}.
\]
Taking $s=r$ and sending $k$ to infinity yields the assertion.
\end{pf}

\section{Proofs of the main results}\label{Proofssection}

\subsection{\texorpdfstring{Proof of Theorem \protect\ref{mainthm}}{Proof of Theorem 2.7}}

We start by making the following simple observation.\eject

\begin{lemma}\label{monotonicitylemma}For any admissible strategy
that is not monotone, there exists another admissible strategy that is
monotone and has strictly lower cost.
\end{lemma}
\begin{pf}We may suppose without loss of generality that $x_0>0$.
Let $x(\cdot)$ be an admissible strategy that is not monotone. Define
$y(t):=x_0+\int_0^t\dot x(t)\Ind{\{\dot x(t)<0\}} \,dt$ and
$\tau:=\inf\{t\ge0 | y(t)=0\}$. Then $y$ is monotone, $\tau\le T$ and
$\widetilde x(t):=y(t\wedge\tau)$ is an admissible strategy with
strictly lower cost than $x(\cdot)$.
\end{pf}

Let $X$ be the superprocess constructed in Proposition
\ref{h-transformexistProp}, and recall the definitions (\ref{vinftydef})
and (\ref{x*}) for $v_\infty$ and $x^*$. In addition to (\ref
{crTdefinition}) we define
%
\begin{equation}
\label{barcrTdefinition} \overline c_{r,T}:=\inf_{t\in[r,T]}
\inf_z\frac{h(t,z)}{\eta
(z)} \quad\mbox{and}\quad C_{r,T}:=
\frac{\overline c_{r,T}}{c_{r,T}}.
\end{equation}
It follows from (\ref{almostharmonic}) that $C_{r,T}$ is strictly
positive for all $r\in[0,T]$ and from (\ref{uniformalmostharmonic})
that $C_{r,T}\to1$ as $r\ua T$. Proposition
\ref{inhomLaplaceestimatesProp}, and hence Theorem \ref{condexpthm},
play a crucial role in proving the following key lemma.

\begin{lemma}\label{finitecostlemma} For $A\in\cA^1_T$ satisfying
(\ref{etaAassumption}), the process
\[
x^*(t)=x_0\exp \biggl(-\int_0^t
\biggl(\frac{v_\infty(s,Z_s)}{\eta
(Z_s)} \biggr)^\beta \,ds \biggr),\qquad 0\le t<T,
\]
is an admissible strategy with finite cost. That is, $x^*(t)\to0$ as
$t\ua T$ and $E_{0,z}[ \int_0^T|\dot x^*(t)|^p\eta(Z_t) \,dt
]<\infty$.
\end{lemma}
\begin{pf} By Lemma \ref{lowerLaplaceestLemma} and (\ref{barcrTdefinition}),
%
\begin{equation}\label{vInftylowerestimate}
v_\infty(r,z)\ge-\log\bP_{r,\delta_z}[ X_T=0 ]\ge
\frac{h(r,z)}{c_{r,T}(T-r)^{1/\beta}}\ge\frac{C_{r,T}\cdot\eta
(z)}{(T-r)^{1/\beta}}.
\end{equation}
We thus get the upper bound
%
\begin{equation}
\label{x*upperbound} x^*(t)\le x_0\exp \biggl(-{C_{r,T}^{\beta}}
\int_r^t\frac1{T-s} \,ds \biggr)=x_0
\biggl(\frac{T-t}{T-r} \biggr)^{C_{r,T}^{\beta}},\qquad r\le t< T.\hspace*{-35pt}
\end{equation}
In particular, we have $x^*(t)\to0$ as $t\ua T$.

Next, recalling the identity $p=\frac{1+\beta}\beta$, estimate
(\ref{x*upperbound}) implies that for a.e. $t\ge r$,
\[
\bigl|\dot x^*(t)\bigr|^p= \biggl(\frac{v_\infty(t,Z_t)}{\eta
(Z_t)} \biggr)^{1+\beta}
x^*(t)^{({\beta+1})/\beta}\le c_1 \biggl(\frac{v_\infty(t,Z_t)}{h(t,Z_t)}(T-t)^{C_{r,T}^{\beta}/\beta}
\biggr)^{1+\beta},
\]
where here and in the sequel $c_i$, $i\in\bN$, denote constants
depending on $r$, $T$, $\beta$, $c_T$, $z$ and $x_0$.
Using (\ref{log-Laplaceestimatewithinfinityeta}) and identity
(\ref{BvsBh}), we obtain that
\begin{eqnarray*}
&& \biggl(\frac{v_\infty(t,Z_t)}{h(t,Z_t)}(T-t)^{C_{r,T}^{\beta
}/\beta} \biggr)^{1+\beta}
\\
&&\qquad\le \biggl( \frac{E_{t,Z_t}[ A[t,T]
]}{h(t,Z_t)}(T-t)^{C_{r,T}^{\beta}/\beta}+c_T(T-t)^{
({C_{r,T}^{\beta}-1})/{\beta}}
\biggr)^{1+\beta}
\\
&&\qquad= \bigl( E^h_{t,Z_t} \bigl[ A_h[t,T]
\bigr](T-t)^{1/\beta
}+c_T \bigr)^{1+\beta}
(T-t)^{({(1+\beta)(C_{r,T}^{\beta
}-1)})/{\beta}}
\\
&&\qquad\le c_2 \bigl(1+E^h_{t,Z_t} \bigl[
A_h[t,T] \bigr]^{1+\beta
}(T-t)^{(1+\beta)/\beta} \bigr)
(T-t)^{({(1+\beta
)(C_{r,T}^{\beta}-1)})/{\beta}}
\\
&&\qquad=c_2(T-t)^{-\delta}+c_2E^h_{t,Z_t}
\bigl[ A_h[t,T] \bigr]^{1+\beta
}(T-t)^{({1+\beta})/\beta-\delta},
\end{eqnarray*}
where
\[
\delta:=-\frac{(1+\beta)(C_{r,T}^{\beta}-1)}{\beta}\ge0.
\]
Now we fix $r$ so that
$\delta<1$,
which is possible since $C_{r,T}\to1$ as $r\ua T$ by (\ref
{uniformalmostharmonic}). Then
%
\begin{eqnarray}
\bigl|\dot x^*(t)\bigr|^p&\le& c_1 \biggl(
\frac{v_\infty
(t,Z_t)}{h(t,Z_t)}(T-t)^{C_{r,t}^\beta/\beta} \biggr)^{1+\beta}\nonumber\\[-8pt]\\[-8pt]
&\le&
c_1c_2(T-t)^{-\delta}+c_1c_3E^h_{t,Z_t}
\bigl[ A_h[t,T] \bigr]^{1+\beta}.\nonumber
\end{eqnarray}

We also need to estimate $|\dot x^*(t)|$ for $0\le t\le r$. To this
end, we will use the trivial bound $x^*(t)\le x_0$ to get as above that
\[
\bigl|\dot x^*(t)\bigr|^p\le c_4 \biggl(\frac{v_\infty
(t,Z_t)}{h(t,Z_t)}
\biggr)^{1+\beta}\le c_5 \bigl(1+E^h_{t,Z_t}
\bigl[ A_h[t,T] \bigr]^{1+\beta} \bigr)
\]
for $0\le t\le r$.

Putting everything together and using $\int_r^T(T-t)^{-\delta}
\,dt<\infty$, the fact that $q=1+\beta$, (\ref{almostharmonic}),
(\ref{BvsBh}), Jensen's inequality, the Markov property of $Z$ and
once again (\ref{almostharmonic}) yields
\begin{eqnarray*}
&& E_{0,z} \biggl[ \int_0^T\bigl|\dot
x^*(t)\bigr|^p\eta(Z_t) \,dt \biggr]
\\
&&\qquad \le c_6c_TT\eta(z)+c_5E_{0,z}
\biggl[ \int_0^r\eta (Z_t)E^h_{t,Z_t}
\bigl[ A_h[t,T] \bigr]^q \,dt \biggr]
\\
&&\qquad\quad{} +c_1c_3E_{0,z} \biggl[ \int
_r^T\eta (Z_t)E^h_{t,Z_t}
\bigl[ A_h[t,T] \bigr]^q \,dt \biggr]
\\
&&\qquad \le c_7+c_8\int_0^TE_{0,z}
\bigl[ \eta(Z_t)E^h_{t,Z_t} \bigl[
A_h[t,T] \bigr]^q \bigr] \,dt
\\
&&\qquad \le c_7+c_8\int_0^TE_{0,z}
\bigl[ \eta (Z_t)h(t,Z_t)^{-q}E_{t,Z_t}
\bigl[ A[t,T]^q \bigr] \bigr] \,dt
\\
&&\qquad \le c_7+c_9\int_0^TE_{0,z}
\bigl[ \eta(Z_t)^{1-q}A[t,T]^q \bigr] \,dt,
\end{eqnarray*}
which is finite due to our assumption (\ref{etaAassumption}).
\end{pf}

\begin{remark}\label{x*boundRemark}
When $\eta\equiv1$ we will also have $C_{r,T}=1$. Hence (\ref
{x*upperbound}) implies the bound (\ref{x*bound}).
\end{remark}

For $p\ge2$, let us introduce the function
%
\begin{equation}
\label{phi} \phi_p(\xi,\zeta):=\xi^p-p
\zeta^{p-1}\xi+(p-1)\zeta^p,\qquad \xi,\zeta\ge0.
\end{equation}
Note that $\phi_p(\xi,\zeta)\ge0$ with equality if and only if $\xi
=\zeta$, because Young's inequality gives
%
\begin{equation}
\label{phipYoungineq} \xi\zeta^{p-1}\le\frac1p\xi^p+
\frac{p-1}p\zeta^p \qquad\mbox{for $\xi$, $\zeta\ge0$.}
\end{equation}

\begin{proposition}\label{costfunctProp}For $A\in\cA^1_T$
satisfying (\ref{etaAassumption}), the following inequality holds
for any monotone admissible strategy $x$ with finite cost:
%
\begin{eqnarray}\label{costidentity}
&& E_{0,z} \biggl[ \int_0^T\bigl|\dot
x(s)\bigr|^p\eta(Z_s) \,ds+\int_{[0,T]}\bigl|x(s)\bigr|^p
A(ds) \biggr]
\nonumber
\\
&&\qquad \ge|x_0|^pv_\infty(0,z)
\\
&&\qquad\quad{} +E_{0,z} \biggl[ \int_0^T
\eta(Z_t)\phi_p \biggl(\bigl|\dot x(t)\bigr|, \bigl|x(t)\bigr| \biggl(
\frac{v_\infty(t,Z_t)}{\eta(Z_t)} \biggr)^{1/(p-1)} \biggr) \,dt \biggr].
\nonumber
\end{eqnarray}
Moreover, there is equality in (\ref{costidentity}) when $x=x^*$.
\end{proposition}

Before proving Proposition \ref{costfunctProp}, let us show how it
implies Theorem \ref{mainthm}.

\begin{pf*}{Proof of Theorem \ref{mainthm}}
Since $\phi_p(\xi,\zeta)\ge0$, Proposition \ref{costfunctProp} yields
that for any monotone admissible strategy $x(\cdot)$ with finite cost,
%
\begin{equation}
\label{costestimateeq} E_{0,z} \biggl[ \int_0^T\bigl|
\dot x(s)\bigr|^p\eta(Z_s) \,ds+\int_{[0,T]}\bigl|x(s)\bigr|^p
A(ds) \biggr]\ge|x_0|^pv_\infty(0,z).
\end{equation}
Since moreover $\phi_p(\xi,\zeta)=0$ if and only if $\xi=\zeta$,
equality holds in (\ref{costestimateeq}) if and only if there is
equality in (\ref{costidentity}) and also
%
\begin{eqnarray}
\label{dotxODE}
\bigl|\dot x(t)\bigr|&=&\bigl|x(t)\bigr| \biggl(\frac{v_\infty(t,Z_t)}{\eta(Z_t)}
\biggr)^{1/(p-1)}\nonumber\\[-8pt]\\[-8pt]
&=&\bigl|x(t)\bigr| \biggl(\frac{v_\infty(t,Z_t)}{\eta(Z_t)} \biggr)^{\beta}
\qquad\mbox{for a.e. $t$, $P_{0,z}$-a.s.}\nonumber
\end{eqnarray}
By the monotonicity of $x$ the latter condition holds if and only if
$x=x^*$. Moreover, $x^*$ is an admissible strategy with finite cost by
Lemma \ref{finitecostlemma} and satisfies equality in~(\ref
{costidentity}).
\end{pf*}

\begin{pf*}{Proof of Proposition \ref{costfunctProp}}
For $k\in\bN$, we introduce the functions
%
\begin{equation}
v_k(r,z)=-\log\bE_{r,\delta_z}\bigl[ e^{-J_A-k\langle1,X_T\rangle} \bigr].
\end{equation}
Then $v_k\nearrow v_\infty$ and $v_k$ uniquely solves
\[
v_k(r,z)=k+E_{r,z}\bigl[ A[r,T] \bigr]-E_{r,z}
\biggl[ \int_r^T v_k(s,Z_s)^{1+\beta}
\frac1{\beta\eta(Z_s)^\beta} \,ds \biggr].
\]

We now fix a monotone admissible strategy $x(\cdot)$ with finite cost.
We may assume without loss of generality that $x_0\ge0$ and hence also
$x(t)\ge0$ for all $t$. We define
%
\begin{equation}
\label{Ckdefinition}\quad C^k_t:=\int_0^t\bigl|
\dot x(s)\bigr|^p\eta(Z_s) \,ds+\int_{[0,t)}x(s)^p
A(ds)+x(t)^pv_k(t,Z_t).
\end{equation}
The first two terms on the right represent the cost accumulated by the
strategy $x$ over the time interval $[0,t)$. The rightmost term
approximates our guess for the minimal cost incurred over the time
interval $[t,T]$ when starting at time $t$ with the remainder $x(t)$.

We next
define $M^k$ as a rightcontinuous version of the martingale
%
\begin{equation}
\label{Mkdefinition} M^k_t:=k+E_{0,z} \biggl[
A[0,T]-\int_0^T v_k(s,Z_s)^{1+\beta}
\frac 1{\beta\eta(Z_s)^\beta} \,ds \Big| \cF_t
\biggr].
\end{equation}
Such a version exists when we replace the filtration $(\cF[0,t])_{t\ge
0}$ by its $P_{0,z}$-augmentation $(\overline\cF[r,t])_{t\ge0}$ so that
the resulting filtered probability space $(\Om,(\overline\cF
[0,\allowbreak t]),P_{0,z})$ satisfies the usual conditions; see Section A.1.1 in
\citet{DynkinBranchingBook}.

By the Markov property of $Z$, we have
%
\begin{equation}
\label{vkdynamics} v_k(t,Z_t)=M^k_t-A[0,t)+
\int_0^t v_k(s,Z_s)^{1+\beta}
\frac1{\beta \eta(Z_s)^\beta} \,ds.
\end{equation}
In particular,
%
\begin{equation}
\label{Vtdef}
V_t:=v_k(t,Z_t)-A\{t\}
\end{equation}
is right continuous, and
%
\begin{equation}
\label{Ckt} C^k_t=\int_0^t\bigl|
\dot x(s)\bigr|^p\eta(Z_s) \,ds+\int_{[0,t]}x(s)^p
A(ds)+x(t)^pV_t
\end{equation}
is rightcontinuous as well.

Our next goal is to investigate the limit of $C^k_t$ as $t\ua T$. To
this end, we define $N$ as a rightcontinuous version of the martingale
\[
N_t=k+E_{0,z}\bigl[ A[0,T] | \cF_t \bigr].
\]
We have
$V_t\le k+E_{t,Z_t}[ A[t,T] ]$ for each $t$ $P_{0,z}$-a.s. and hence,
by right continuity, $V_t\le N_t$ for all $t$ $P_{0,z}$-a.s.
The martingale convergence theorem implies that $\sup_{t\le
T}N_t<\infty$ $P_{0,z}$-a.s. Hence also $\sup_{t\le T}V_t<\infty$
$P_{0,z}$-a.s. It thus follows from $x(t)\to0$ that $x(t)^pV_t\to0$
$P_{0,z}$-a.s. as $t\ua T$.
Moreover, we have $\int_{[0,T)}x(s) A(ds)=\int_{[0,T]}x(s) A(ds)$
because $x(T)=0$. Therefore,
%
\begin{equation}
\label{Ctconvergence} C^k_t \lra C^k_T:=
\int_0^T\bigl|\dot x(s)\bigr|^p
\eta(Z_s) \,ds+\int_{[0,T]}x(s)^p A(ds),
\end{equation}
$P_{0,z}$-a.s. as $t\ua T$.

Next, applying It\^o's formula to (\ref{Ckt}) and using (\ref
{vkdynamics}) and (\ref{Vtdef}) yields
\begin{eqnarray*}
dC^k_t &=&\bigl|\dot x(t)\bigr|^p
\eta(Z_t) \,dt+x(t)^p A(dt)+px(t)^{p-1}\dot
x(t)V_t \,dt+x(t)^p \,dV_t
\\
&=& \biggl(\bigl|\dot x(t)\bigr|^p\eta(Z_t)+px(t)^{p-1}
\dot x(t)V_t+ x(t)^pV_t^{1+\beta}\frac1{
\beta\eta(Z_t)^\beta} \biggr) \,dt\\
&&{}+x(t)^p
\,dM^k_t
\\
&=&\eta(Z_t)\phi_p \biggl(\bigl|\dot x(t)\bigr|, x(t) \biggl(
\frac
{v_k(t,Z_t)}{\eta(Z_t)} \biggr)^{1/(p-1)} \biggr) \,dt+x(t)^p
\,dM^k_t,
\end{eqnarray*}
where, in the last step, we have used the relation $1+\beta=p/(p-1)$
and the fact that $V_t=v_k(t,Z_t)$ for a.e. $t$.

Using (\ref{Ctconvergence}), we obtain in the limit $t\ua T$ that
$P_{0,z}$-a.s.
\begin{eqnarray*}
&& \int_0^T\bigl|\dot x(s)\bigr|^p \,ds+\int
_{[0,T]}x(s)^p A(ds)-x_0^pv_k(0,z)
\\
&&\qquad=C^k_T -C^k_0
\\
&&\qquad=\int_0^T\eta(Z_t)
\phi_p \biggl(\bigl|\dot x(t)\bigr|, x(t) \biggl(\frac
{v_k(t,Z_t)}{\eta(Z_t)}
\biggr)^{1/(p-1)} \biggr) \,dt+ \int_0^Tx(t)^p
\,dM^k_t.
\end{eqnarray*}
Taking expectations yields
%
\begin{eqnarray}\label{phipvkidentity}
&& E_{0,z} \biggl[ \int_0^T\bigl|\dot
x(s)\bigr|^p\eta(Z_s) \,ds+\int_{[0,T]}x(s)^p
A(ds) \biggr]
\nonumber\\
&&\qquad =x_0^pv_k(0,z)\\
&&\qquad\quad{}+E_{0,z} \biggl[
\int_0^T\eta(Z_t)\phi_p
\biggl(\bigl|\dot x(t)\bigr|, x(t) \biggl(\frac{v_k(t,Z_t)}{\eta(Z_t)} \biggr)^{1/(p-1)} \biggr)
\,dt \biggr].\nonumber
\end{eqnarray}
Fatou's lemma and the facts that $\phi_p\ge0$ and $v_k\ua v_\infty$
thus yield the first part of the assertion when
passing to the limit $k\ua\infty$ on the right-hand side of this identity.

Now we show that we may retain equality in (\ref{phipvkidentity})
when passing to the limit $k\ua\infty$ and taking $x=x^*$. By (\ref
{dotxODE}), (\ref{phipYoungineq}) and dominated convergence, this
will hold when
\[
E_{0,z} \biggl[ \int_0^T
\eta(Z_t) \biggl( x^*(t) \biggl(\frac{v_\infty
(t,Z_t)}{\eta(Z_t)}
\biggr)^{1/(p-1)} \biggr)^p \,dt \biggr]=E_{0,z} \biggl[
\int_0^T\bigl|\dot x^*(t)\bigr|^p
\eta(Z_t) \,dt \biggr]
\]
is finite. But this is true by Lemma \ref{finitecostlemma}, and so
the second part of the assertion is proved.
\end{pf*}

\subsection{\texorpdfstring{Proof of Theorem \protect\ref{relaxedthm}}{Proof of Theorem 2.8}}

\begin{lemma}\label{vvarrhoconvergencelemma}
If $\varrho$ is as in Theorem \ref{relaxedthm}, then
$
v_\varrho(t,Z_t)\to\varrho(Z_T)$ in $P_{0,z}$-probability as $t\ua T$.
\end{lemma}
\begin{pf}
First, using the Markov property of $Z$,
%
\begin{eqnarray}
\label{lemma1steq} v_\varrho(t,Z_t)&\le& E_{t,Z_t}
\bigl[ A[t,T]+\varrho(Z_T) \bigr]
\nonumber\\[-8pt]\\[-8pt]
&=&E_{0,z} \bigl[ A[t,T] | \cF[0,t] \bigr]+E_{t,Z_t}\bigl[
\varrho (Z_T) \bigr].\nonumber
\end{eqnarray}
By (\ref{varrhocond}), the second term on the right converges
$P_{0,z}$-a.s. to $\varrho(Z_T)$. As for the first term on the right,
we first note that
$A[t,T]\to0$ in $L^1(P_{0,z})$ due to our assumption $A\{T\}=0$
$P_{0,z}$-a.s., the fact that $E_{0,z}[ A[0,T] ]<\infty$, and
dominated convergence. Therefore, $E_{0,z}[ A[t,T] | \cF[0,t] ]\to
0$ in $P_{0,z}$-probability, and so the entire right-hand side of
(\ref{lemma1steq}) converges to $\varrho(Z_T)$ in $P_{0,z}$-probability.

We also need a lower bound on $v_\varrho(t,Z_t)$. Using the results
and notation from Section \ref{h-trafosection}, we have
%
\begin{eqnarray}
\label{varrholemma2ndeq} v_\varrho(t,Z_t)&=&-\log
\bE_{t,\delta_{Z_t}}\bigl[ e^{-J_A-\langle
\varrho,X_T\rangle} \bigr]\ge-\log\bE_{t,\delta_{Z_t}}\bigl[
e^{-\langle
\varrho,X_T\rangle} \bigr]
\nonumber\\[-8pt]\\[-8pt]
&=&-\log\bE^h_{t,h(t,Z_t)\delta_{Z_t}}\bigl[ e^{-\langle\varrho/\eta,X^h_T\rangle}
\bigr]=h(t,Z_t)v^h(t,Z_t),\nonumber
\end{eqnarray}
where $v^h$ solves
\[
v^h(r,y)=E^h_{r,y} \biggl[
\frac{\varrho(Z_T)}{\eta(Z_T)} \biggr]-E^h_{r,y} \biggl[ \int
_r^T\frac1\beta \biggl(\frac{h(t,Z_t)}{\eta
(Z_t)}
\biggr)^\beta v^h(t,Z_t)^{1+\beta} \,dt
\biggr].
\]
By (\ref{almostharmonic}), we have $h(t,Z_t)/\eta(Z_t)\le c_T$.
Let $\widetilde X=(\widetilde X_t,\cG(I),\widetilde\bP_{r,\mu})$
be the superprocess
with one-particle motion $Z^h$, $K(ds)=\frac{c_T^\beta}{\beta} \,ds$,
and $\psi(\xi)=\xi^{1+\beta}$. Due to the fact that $E^h_{r,y}[
\frac{\varrho(Z_T)}{\eta(Z_T)} ]=\frac1{h(r,y)}E_{r,y}[ \varrho
(Z_T) ]<\infty$ for all $(r,y)$ by (\ref{almostharmonic}) and the
first condition in (\ref{varrhocond}), we may apply Proposition
\ref{MaxprinzProp}, which then yields that
\[
v^h(r,y)\ge-\log\widetilde\bE_{r,\delta_y}\bigl[
e^{-\langle\varrho
/\eta,\widetilde X_T\rangle} \bigr].
\]
Since $\widetilde X$ has a homogeneous branching mechanism, we may
apply the
lower bound from (\ref{Laplaceineq}) to get
\begin{eqnarray*}
-\log\widetilde\bE_{r,\delta_y}\bigl[ e^{-\langle\varrho/\eta,\widetilde
X_T\rangle} \bigr]
&\ge&E_{r,y}^h \biggl[ \frac{{\varrho
(Z_T)}/{\eta(Z_T)}}{(1+c_T^\beta(T-r) ({\varrho(Z_T)}/{\eta
(Z_T)})^\beta )^{1/\beta}} \biggr]
\\
&=&\frac1{h(r,y)}E_{r,y} \biggl[ \frac{\varrho(Z_T)}{(1+c_T^\beta
(T-r) ({\varrho(Z_T)}/{\eta(Z_T)})^\beta )^{1/\beta}} \biggr]
\\
&\ge&\frac{E_{r,y}[ \varrho(Z_T) ]}{h(r,y)(1+c_T^\beta c_\varrho
^\beta(T-r))^{1/\beta}},
\end{eqnarray*}
where we have used the first condition from (\ref{varrhocond}) in the
third step. Combining the preceding inequality with (\ref
{varrholemma2ndeq}) yields
\[
\liminf_{t\ua T}v_\varrho(t,Z_t)\ge\liminf
_{t\ua T}\frac
{E_{t,Z_t}[ \varrho(Z_T) ]}{(1+c_T^\beta c_\varrho^\beta
(T-t))^{1/\beta}}=\varrho(Z_T),\qquad
\mbox{$P_{0,z}$-a.s.}
\]
and hence the assertion.
\end{pf}

\begin{pf*}{Proof of Theorem \ref{relaxedthm}}
Just as in Lemma \ref{monotonicitylemma} one first notes that we can
restrict our attention to monotone relaxed strategies $x(\cdot)$ with
$x(T)\ge0$ for $x_0\ge0$ and $x(T)\le0$ for $x_0\le0$.

We may assume $x_0\ge0$ without loss of generality. For a given
monotone relaxed strategy $x(\cdot)$ with $x(T)\ge0$ we define
\[
C_t:=\int_0^t\bigl|\dot
x(s)\bigr|^p\eta(Z_s) \,ds+\int_{[0,t)}x(s)^p
A(ds)+x(t)^pv_\varrho(t,Z_t).
\]
It follows from our assumptions $A\{T\}=0$ and Lemma \ref
{vvarrhoconvergencelemma} that in $P_{0,z}$-probability
\[
C_t\lra C_T=\int_0^T\bigl|
\dot x(s)\bigr|^p\eta(Z_s) \,ds+\int_{[0,T]}x(s)^p
A(ds)+x(t)^p\varrho(Z_T).
\]

The function $v_\varrho$ solves
\[
v_\varrho(r,z)=E_{r,z} \bigl[ A[r,T]+\varrho(Z_T)
\bigr]-E_{r,z} \biggl[ \int_r^T
v_\varrho(s,Z_s)^{1+\beta}\frac1{\beta
\eta(Z_s)^\beta } \,ds \biggr].
\]
Thus, arguing as in the proof of Proposition \ref{costfunctProp}, we
find that
\begin{eqnarray*}
&& E_{0,z} \biggl[ \int_0^T\bigl|\dot
x(s)\bigr|^p\eta (Z_s) \,ds+\int_{[0,T]}x(s)^p
A(ds)+\varrho(Z_T)x(T)^p \biggr]
\\
&&\qquad=x_0^pv_\varrho(0,z)+E_{0,z} \biggl[
\int_0^T\eta(Z_t)
\phi_p \biggl(\bigl|\dot x(t)\bigr|, x(t) \biggl(\frac{v_\varrho(t,Z_t)}{\eta(Z_t)}
\biggr)^{1/(p-1)} \biggr) \,dt \biggr].
\end{eqnarray*}
In contrast to the proof of Proposition \ref{costfunctProp}, note
that here is no need for a limiting procedure since, unlike $v_\infty
$, the function $v_\varrho$ has no singularity.

Now we can proceed as in the proof of Theorem \ref{mainthm} to get
the assertion.
\end{pf*}

\section*{Acknowledgment}

The author wishes to thank Aur\'elien Alfonsi and an anonymous referee
for helpful comments on a previous draft of the paper.



\printaddresses

\end{document}